 \input amstex
  \documentstyle{amsppt}
  \NoBlackBoxes
  \pagewidth{137mm}
  \pageheight{198mm}
 \loadbold
 \magnification=1200
 \baselineskip=15 pt
 \TagsOnRight
 %
 \define\CC{\Cal C}
 \define\UU{\Cal U}
 \define\FF{\Cal F}
 \define\MM{\Cal M}
 \define\RR{\Cal R}
 \define\BC{\Bbb C}
 \define\BZ{\Bbb Z}
 \define\BI{\Bbb I}
 \define\vphi{\varphi}
 
 \define\veps{\varepsilon}
 \define\sg{\sigma}
 \define\gtg{\goth{g}}
 \define\gk{\goth{k}}
 \define\gh{\goth{h}}
 \define\gb{\goth{b}}
 \define\gn{\goth{n}}
 \define\gX{\goth{X}}
 \define\gsl{\goth{sl}}
 \define\gR{\goth{R}}
 \define\tx{\tilde\xi}
 \define\tp{\tilde\varphi}
 \define\te{\tilde\varepsilon}
 \define\td{\tilde\Delta}
 \define\z{\bar{z}}
 \define\lo{_{(1)}}
 \define\lt{_{(2)}}
 \define\loo{_{(1)(1)}}
 \define\lot{_{(1)(2)}}
 \define\lto{_{(2)(1)}}
 \define\ltt{_{(2)(2)}}
 \define\Zp{\BZ_{+}}
 \define\Uq{\Cal U_q}
 \define\Lp{L^{+}}
 \define\Lm{L^{-}}
 \define\Lpm{L^{\pm}}
 
 \define\Z{Z^{\ast}}
 
 \define\spn{\operatorname{span}}
 
 \define\diag{\operatorname{diag}}
 \define\id{\operatorname{id}}
 \define\ad{\operatorname{ad}}
 \define\gr{\operatorname{gr}}
 \define\Funq{\operatorname{Fun}_{q}}
 \define\idr{\langle\Cal R\rangle}
 
 %
 \define\JMPI{1}
 \define\Grass{2}
 \define\JSII{3}
 \define\JAlg{4}
 \define\JMPII{5}
 \define\Bulgariens{6}
 \define\JurcoSch{7}
 \define\ParshalWang{8}
 \define\TauberTaft{9}
 \define\BiedenharnL{10}
 \define\Dabrowetal{11}
 \define\Drinfeld{12} 
 \define\Jimbo{13}    
 \define\Podlesz{14}
 \define\Sheu{15}
 \define\Bayenetal{16}
 \define\Zuminoetal{17}
 \define\Rosso{18}
 \define\Lusztig{19}
 \define\FRT{20}
 \define\JimboR{21} 
 \define\DrinfeldP{22} 
 \define\CJP{23} 
 \define\Joseph{24}
 %
 \topmatter
 \title A construction of representations and quantum homogeneous spaces
 \endtitle
 \author Pavel \v S\v tov\'\i\v cek
 \endauthor
 \affil
 Department of Mathematics\\
 Faculty of Nuclear Science, CTU\\
 Trojanova 13, 120 00 Prague, Czech Republic\\
 stovicek\@kmdec.fjfi.cvut.cz
 \endaffil
 \abstract
A simplified construction of representations  is presented for the
quantized enveloping algebra $\Uq(\gtg)$, with $\gtg$ being a simple
complex Lie algebra belonging to one of
the four principal series $A_\ell$, $B_\ell$, $C_\ell$ or $D_\ell$.
The carrier representation space is the quantized algebra of polynomials
in antiholomorphic coordinate functions on the big cell of a coadjoint
orbit of $K$ where $K$ is the compact simple Lie group with the Lie
algebra $\gk$ -- the compact form of $\gtg$.
 \endabstract
 \endtopmatter
 \document
 \bigpagebreak
 \flushpar {\bf 1. Motivation} 
 \medpagebreak

Let $G$ be a simple and simply connected Lie group belonging to one of
the four principal series $A_\ell$, $B_\ell$, $C_\ell$ or $D_\ell$, and
$K\subset G$ its compact form. The symbols $\gtg$ and $\gk$ designate
the corresponding Lie algebras and $\ell$ is equal to the rank. The
primary motivation was to find a quantum version of the construction of
representations for the group $K$ via the method of orbits due to
Kirillov and Kostant, with the result being expressed explicitly in
terms of local holomorphic (or antiholomorphic) coordinates on the
coadjoint orbit. But first let us consider briefly the classical case.

Each $G$ is a complex matrix group. The tautological
(defining) representation $T$ is frequently called the vector
representation. Let us denote by $N$ its dimension (i.e.,
$G\subset SL(N,\BC)$).
Every coadjoint orbit $X=K_0\backslash K=P_0\backslash G$ is a
homogeneous space for both $K$ and $G$ and so it is a compact complex
manifold. We shall restrict ourselves to the generic orbits of the top
dimension $\dim_\BC\,X=(\dim_\BC\,\gtg-\ell)/2$. The local holomorphic
coordinates on $X$ are introduced with the aid of Gauss decomposition.
The factor mapping $G\to X=P_0\backslash G$ sends an element
$g\in G$ to a point belonging to the so called big cell if and only
if there exists a decomposition
$g=g_{(-)}Z$
where $g_{(-)}$ is a lower triangular matrix and $Z$ is upper triangular
with units on the diagonal,
$$
Z=\pmatrix
1 & z_{12} & \hdots & z_{1N} \\
0 & 1 & \hdots & z_{2N} \\
\hdotsfor 4 \\
\hdotsfor 4 \\
0 & 0 & \hdots & 1
\endpmatrix \,.
\tag 1
$$
In the case of the series $B_\ell$, $C_\ell$ and $D_\ell$ the subgroup
$G\subset SL(N,\BC)$ is determined by the equation
$C_0 g^t C_0^{-1} = g^{-1}$
where $C_0$ is an appropriate real $N\times N$ matrix
($(C_0)_{jk}=\pm\delta_{j+k,N+1}$). Consequently the matrix
$Z$ must obey a similar condition,
$$
C_0 Z^t C_0^{-1} = Z^{-1}.
\tag 2
$$
In fact, the equality (2) reduces the number of independent coordinates
$z_{jk}$ living on the big cell to the correct value
$(\dim_\BC\,\gtg-\ell)/2$.

As usual, one makes use of the fact that
irreducible unitary representations of $K$ are in one-to-one
correspondence with finite-dimensional irreducible representations
of the Lie algebra $\gtg$ over $\BC$.
The right action of $G$ on  the big cell
can be described, too, with the help
of Gauss decomposition. For $g\in G$ and $Z$ as above let us
decompose, if possible,
$Zg=(Zg)_{(-)}(Zg)_{(+)}$
where again $(Zg)_{(-)}$ is lower triangular and $(Zg)_{(+)}$ is upper
triangular with units on the diagonal. The right action reads
$$
Z\cdot g := (Zg)_{(+)}.
\tag 3
$$
Naturally, the right hand side of (3) is not well defined for all $g$
and $Z$, i.e., it has singularities. The reason is simple -- the
coordinates $z_{jk}$ are local while the action itself is global.
However by differentiating the equality (3)
one gets a well defined infinitesimal action
$\xi:\gtg\to\gX_\BC(X)$. Each element $x\in\gtg$ is represented by
a complex vector field $\xi_x$, with $\xi_x$ depending on $x$ linearly
over $\BC$, and it holds
$[\xi_x,\xi_y]=\xi_{[x,y]}$.

The infinitesimal action $\xi$ doesn't lead directly to a
finite-dimensional irreducible representation of $\gtg$. The result of
the method of orbits, when looking at it through the local coordinates,
is a correction achieved
by adding to $\xi_x$ a holomorphic function $\vphi_x$
defined on the big cell (in fact, a polynomial in the coordinates
$z_{jk}$) and depending on $x$ linearly. Thus elements from $\gtg$ are
represented by first order differential operators:
$x\mapsto\xi_x+\vphi_x $.
The Lie bracket is preserved provided the function $\vphi_x$ fulfills
the condition
$$
\xi_x\cdot\vphi_y-\xi_y\cdot\vphi_x=\vphi_{[x,y]}.
\tag 4
$$
The carrier vector space of the representation is built up from
holomorphic functions on the big cell (in fact, from polynomials in
$z_{jk}$), with the unit function as a cyclic vector.

A comparatively simple construction presented in this letter attempts to
generalize the experience accumulated throughout the series of papers
\cite{\JMPI, \Grass, \JSII, \JAlg, \JMPII} and to simplify the
procedures applied therein as much as possible. One should mention also
the papers \cite{\Bulgariens, \JurcoSch} which were prior to the paper
\cite{\JMPII}. Even more, this construction makes it possible to deal
with the general case, including the orthogonal and symplectic groups
(the series $B_\ell$, $C_\ell$ and $D_\ell$). Up to now, only some
particular cases were treated in \cite{\JAlg}.
The main idea of the construction consists in finding a quantum analog
to the function $\vphi$, including the compatibility condition (4).

Of course, apart of this local point of view there are known
also other algebraic constructions, particularly global ones making use
of the idea of induced representations and with the carrier
representation space being formed by holomorphic sections in quantized
line bundles over $X$. Let us mention just a few papers dealing with
this subject: \cite{\ParshalWang, \TauberTaft, \BiedenharnL,
\Dabrowetal}.

 \bigpagebreak
 \flushpar {\bf 2. Construction}
 \medpagebreak

Let us introduce the initial data. Assume we are given a
bialgebra $\UU$ with the counit denoted by $\veps$ and the
comultiplication denoted by $\Delta$ (the antipode will not be used and
so $\UU$ need not be a Hopf algebra) and a unital algebra
$\CC$. Moreover, $\CC$ is supposed to be a left $\UU$-module with the
action denoted by $\xi$,
$$
\UU\otimes\CC\ni x\otimes f\mapsto\xi_x\cdot f\in\CC,
\tag 5
$$
and fulfilling two conditions:
$$
\align
& \xi_x\cdot 1=\veps(x)\,1,\quad \forall x\in\UU,
\tag 6\\
& \xi_x\cdot(fg)=(\xi_{x\lo}\cdot f)(\xi_{x\lt}\cdot g),\quad
\forall x\in\UU,\ \forall f,g\in\CC.
\tag 7
\endalign
$$
If convenient we shall write $\xi(x)\cdot f$ instead of
$\xi_x\cdot f$. The second condition (7) is nothing but Leibniz
rule. Here and everywhere in what follows we use Sweedler's notation:
$\Delta x=x\lo\otimes x\lt$. Note that the coassociativity of
$\Delta$ can be expressed in this formalism as
$$
x\loo\otimes x\lot\otimes x\lt=x\lo\otimes x\lto\otimes x\ltt.
\tag 8
$$

\proclaim{Proposition 1}
Suppose that a linear mapping $\vphi:\UU\to\CC$ satisfies
$$
\align
& \vphi(1)=1,
\tag 9\\
& \vphi(xy)=\bigl(\xi_{x\lo}\cdot\vphi(y)\bigr)\vphi(x\lt),
\quad\forall x,y\in\UU.
\tag 10
\endalign
$$
Then the prescription
$$
x\cdot f:=(\xi_{x\lo}\cdot f)\,\vphi(x\lt),
\quad\forall x\in\UU,\ \forall f\in\CC,
\tag 11
$$
defines a left $\UU$-module structure on $\CC$ and it holds
$$
x\cdot(fg)=(\xi_{x\lo}\cdot f)(x\lt\cdot g),\quad
\forall x\in\UU,\ \forall f,g\in\CC.
\tag 12
$$
Particularly,
$$
\vphi(x)=x\cdot 1,\quad\forall x\in\UU.
\tag 13
$$

Conversely, suppose that
$\UU\otimes\CC\to\CC:x\otimes f\mapsto x\cdot f$
is a left $\UU$-module structure on $\CC$ such that the rule (12) is
satisfied. Then the linear mapping $\vphi:\UU\to\CC$ defined by the
equality (13) fulfills (9) and (10), and consequently the prescription
(11) holds true.
\endproclaim

\remark{Remarks}
The condition (10) generalizes (4). Moreover, (9) "almost" follows from
(10). More precisely, set $x=1$ in (10) to get the equality
$\vphi(y)=\vphi(y)\vphi(1)$, $\forall y\in\UU$.
So (9) is a consequence of (10) as soon as there exists at least one
element $y\in\UU$ such that $\vphi(y)$ is not a left divisor of zero.

The property (12) can be regarded as a generalized Leibniz rule.
\endremark

\demo{Proof}
Let us consider only the first part of the proposition. All
verifications are quite straightforward. We have
$$
1\cdot f=(\xi_1\cdot f)\,\vphi(1)=f
\tag 14
$$
and
$$
\aligned
x\cdot(y\cdot f) &= \left({\xi_{x\lo}\cdot\bigl(
(\xi_{y\lo}\cdot f)\,\vphi(y\lt)\bigr)}\right)\,
   \vphi(x\lt) \\
&= \bigl(\xi_{x\loo}\cdot(\xi_{y\lo}\cdot f)\bigr)
   \bigl(\xi_{x\lot}\cdot\vphi(y\lt)\bigr)\,\vphi(x\lt) \\
&= (\xi_{x\lo y\lo}\cdot f)
   \bigl(\xi_{x\lto}\cdot\vphi(y\lt)\bigr)\,\vphi(x\ltt) \\
&= (\xi_{x\lo y\lo}\cdot f)\,\vphi(x\lt y\lt) \\
&= (xy)\cdot f.
\endaligned
\tag 15
$$
Furthermore,
$$
\aligned
x\cdot(fg) &= \bigl(\xi_{x\lo}\cdot(fg)\bigr)\,\vphi(x\lt) \\
&= (\xi_{x\loo}\cdot f)(\xi_{x\lot}\cdot g)\,\vphi(x\lt) \\
&= (\xi_{x\lo}\cdot f)(\xi_{x\lto}\cdot g)\,\vphi(x\ltt) \\
&= (\xi_{x\lo}\cdot f)(x\lt\cdot g).
\endaligned
\tag 16
$$
Finally,
$$
\align
x\cdot 1 &= (\xi_{x\lo}\cdot 1)\,\vphi(x\lt)
= \vphi(\veps(x\lo)\,x\lt) = \vphi(x).
\tag 17\\
& \tag"{$\qed$}"
\endalign
$$
\enddemo

Here is a trivial example. Suppose that $\chi:\UU\to\BC$ is an algebra
homomorphism (a character) and set $\vphi(x)=\chi(x)\,1$,
$\forall x\in\UU$. Then $\vphi$ fulfills both (9) and (10).
Particularly, for $\chi=\veps$ we recover the original $\UU$-module
structure, i.e., $x\cdot f=\xi_x\cdot f$.

To deal with the function $\vphi$ let us suppose, as usual, that $\UU$
is generated as an algebra by a set of generators $\MM\subset\UU$. Let
$\FF$ be the free algebra generated by $\MM$. Thus $\UU$ is identified
with a quotient $\FF/\idr$ where $\idr$ is the ideal generated by a
set of defining relations $\RR\subset\FF$. Let $\pi$ be the factor
morphism, $\pi:\FF\to\UU$.
Using $\pi$ one can pull back various structures from $\UU$ to $\FF$.
Particularly we set
$$
\gather
\te:=\veps\circ\pi,
\tag 18 \\
\tx_x\cdot f:=\xi_{\pi(x)}\cdot f,\quad\forall x\in\FF,\
\forall f\in\CC .
\tag 19
\endgather
$$
Clearly, $\tx(\idr)=0$.

In addition we assume that the set of generators $\MM\subset\UU$ behaves
well with respect to the comultiplication. More precisely, we impose the
condition
$$
\Delta(\MM)\subset\spn_\BC(\MM_1\otimes\MM_1)\quad
\text{where}\quad\MM_1:=\MM\cup\{1\}.
\tag 20
$$
This property is also fulfilled quite frequently. Then it is natural to
define a comultiplication $\td$ on $\FF$ by the equality
$$
\td(x_1\dots x_n):=\Delta(x_1)\dots\Delta(x_n),\quad
x_i\in\MM.
\tag 21
$$
As $\UU$ is a bialgebra $\RR$ must satisfy
$$
\td(\RR)\subset\idr\otimes\FF+\FF\otimes\idr.
\tag 22
$$
In other words $\idr$ is, at the same time, a coideal.

It is not difficult to check that $\FF$ becomes this way a bialgebra and
that the triple $(\FF,\tx,\CC)$ fulfills the original conditions (6)
and (7), just replacing $\UU$ with $\FF$ and $\xi$ with $\tx$.
Hence Proposition 1 can be applied to $\FF$ as well. But in the case of
the free algebra it is rather easy to describe all linear mappings
$\tp:\FF\to\CC$ with the desired properties.

\proclaim{Lemma 2}
Let the symbols $\FF$ and $\MM$ have the same meaning as above. Then to
any mapping $\vphi:\MM\to\CC$ there exists a unique linear extension
$\tp:\FF\to\CC$ such that $\tp(1)=1$  and the property
$$
\tp(xy)=\bigl(\tx_{x\lo}\cdot\tp(y)\bigr)\tp(x\lt),
\tag 23
$$
is satisfied for all $x,y\in\FF$.
\endproclaim

\demo{Proof}
The algebra $\FF$ is naturally graded,
$$
\FF={\sum_{n\in\Zp}}\sp{\!\!\bigoplus}\,\,\FF^{(n)}\quad
\text{where}\quad\FF^{(n)}:=\spn_\BC(\MM^n).
\tag 24
$$
The mapping $\tp$ is prescribed on $\FF^{(0)}\oplus\FF^{(1)}$.  There is
no contradiction here since (23) is automatically satisfied as son as
$x=1$ and $y$ is arbitrary (obvious) or $x$ is arbitrary and $y=1$ (the
same verification as in (17)). One can proceed by induction in $n$ in
order to extend the definition of $\tp$ to all monomials
$x_1\dots x_n\in\MM^n$ for all $n\in\Zp$. The induction step
$n\to n+1$ is dictated by the rule (23). Thus we set
$$
\tp(xy):= \bigl(\tx_{x\lo}\cdot\tp(y)\bigr)\tp(x\lt),
\quad\forall x\in\MM,\ \forall y\in\MM^n.
\tag 25
$$
Note that the definition (25) makes good sense owing to the property
(20). This means that there is a unique way how to extend $\tp$ from
$\FF^{(0)}\oplus\FF^{(1)}$ to the whole algebra $\FF$ so that $\tp$ is
linear and the rule (23) holds true for all $x\in\MM$ and all $y\in\FF$.
To finish the proof we have to show that this rule is actually satisfied
for all $x\in\MM^n$ and all $y\in\FF$, with $n\in\Zp$ being arbitrary.
We shall again proceed by induction in $n$. In order to carry out the
induction step $n\to n+1$ let us suppose that $x\in\MM$, $z\in\MM^n$,
and that it holds
$\tp(zy) = \bigl(\tx_{z\lo}\cdot\tp(y)\bigr)\tp(z\lt)$,
$\forall y\in\FF$. Then
$$
\aligned
\tp(xzy) &= \left({\tx_{x\lo}\cdot\bigl(
(\tx_{z\lo}\cdot\tp(y))\,\tp(z\lt)\bigr)}\right)\,
   \vphi(x\lt) \\
&= \bigl(\tx_{x\loo}\cdot(\tx_{z\lo}\cdot\tp(y) )\bigr)
   \bigl(\tx_{x\lot}\cdot\tp(z\lt)\bigr)\,\tp(x\lt) \\
&= (\tx_{x\lo z\lo}\cdot \tp(y) )
   \bigl(\tx_{x\lto}\cdot\tp(z\lt)\bigr)\,\tp(x\ltt) \\
&= (\tx_{x\lo z\lo}\cdot \tp(y) )\,\tp(x\lt z\lt).
\endaligned
\tag 26
$$
Thus we have verified that (23) holds also true with $x$ being replaced by
$xz$.
\qed
\enddemo

The final step of the construction is to decide under which conditions
the mapping $\tp$ admits factorization from $\FF$ to $\UU=\FF/\idr$.

\proclaim{Proposition 3}
Suppose there is given a mapping $\vphi:\MM\to\CC$ and let
$\tp$ be its extension to $\FF$ as described in Lemma 2. If
$$
(\pi\otimes\tp)\circ\td(\RR)=0
\tag 27
$$
then $\tp(\idr)=0$ and so there exists a unique linear mapping
$\vphi':\UU\to\CC$ such that $\tp=\vphi'\circ\pi$. Moreover, $\vphi'$
satisfies the conditions (9) and (10).

The same conclusions hold true provided $\RR$ fulfills a stronger condition
than (22), namely
$$
\td(\RR)\subset\idr\otimes\FF+\FF\otimes\FF\RR,
\tag 28
$$
while $\tp$ satisfies a weaker condition
$$
\tp(\RR) = 0.
\tag 29
$$
\endproclaim

\remark{Remark}
Notice that $\tp'|\MM=\vphi|\MM$ and so we have the right to
suppress the dash in the notation and this is what we shall do from
now on.
\endremark

\demo{Proof}
We have to show that
$\tp(xyz)=0$ holds true for all $x,z\in\FF$ and all $y\in\RR$.
According to Proposition 1 there is a unique left $\FF$-module structure
on $\CC$ associated with $\tp$ and we have $\tp(u)=u\cdot 1$,
$\forall u\in\FF$. Suppose that $x,z\in\FF$ and $y\in\RR$. The assumption
(27) then implies
$$
(xyz)\cdot1=x\cdot\left(\bigl(\tx_{y\lo}\cdot\tp(z)\bigr)\,
\tp(y\lt)\right)=0.
\tag 30
$$
This verifies the existence of a linear mapping
$\vphi:\UU\to\CC$ as claimed.

More generally, we find that
$x\cdot f=(\tx_{x\lo}\cdot f)\,\tp(x\lt)=0$
holds true for all $x\in\idr$ and all $f\in\CC$. Consequently the left
action $\FF\otimes\CC\to\CC$ can be factorized from $\FF$ to
$\UU=\FF/\idr$ and the obtained action $\UU\otimes\CC\to\CC$ is
associated with the found mapping $\vphi:\UU\to\CC$. Proposition 1
(the second part) then concludes the proof.

As for the second part one can employ the property (28) to show that
the condition (27) follows from (29).
\qed
 \enddemo

 \bigpagebreak
 \flushpar {\bf 3. Examples}
 \medpagebreak

Now we are going to apply the construction described in the previous
section to the case when
$\UU=\Uq(\gtg)$ is the quantized enveloping algebra
in the sense of Drinfeld \cite{\Drinfeld} and Jimbo \cite{\Jimbo}, and
$\CC$ is the quantized big cell of the orbit
$X=K_0\backslash K=P_0\backslash G$ considered as a complex manifold
\cite{\JSII}. In what follows we prefer the antiholomorphic coordinates
to the holomorphic ones. We assume that the deformation parameter $q>0$,
$q\neq1$, and we set $[x]:=(q^x-q^{-x})/(q-q^{-1})$ for any $x\in\BC$.

 \medpagebreak
 \flushpar {\it 3.1 The particular case $\gtg=\gsl(2,\BC)$}
 \smallpagebreak

First we wish to treat separately the simplest particular case when
$\gtg=\gsl(2,\BC)$ and $X$ is just the Riemannian sphere
(one-dimensional complex projective space). This basic example of a
quantum homogeneous space was analyzed by Podlesz \cite{\Podlesz}.
Afterwards it was reconsidered several times from various points of view.
For example, it fits well the general scheme of deformation quantization
\cite{\Sheu} as proposed in \cite{\Bayenetal}. In our approach we
prefer the local description in terms of coordinates on the big cell
denoted by $z,\z$ \cite{\JMPI}. This point of view was further developed
in \cite{\Zuminoetal}. However when restricting ourselves to the
antiholomorphic part we are left with the algebra of polynomials
$\CC=\BC [\z]$, the same one as in the classical case.

We choose the standard set of generators
$\MM=\{q^{H/2},q^{-H/2},X^+,X^-\}$. The defining relations are also the
usual ones:
$$
\gathered
q^{H/2}q^{-H/2}-1=q^{-H/2}q^{H/2}-1=
q^{H/2}X^\pm-q^{\pm1}X^\pm q^{H/2}=0, \\
[X^+,X^-]-(q-q^{-1})^{-1}(q^H-q^{-H})=0
\endgathered
\tag 31
$$
(of course, $q^{\pm H}\equiv(q^{\pm H/2})^2$). Let us recall, too, the
formulae for the comultiplication and the counit:
$$
\gathered
\Delta(q^{\pm H/2}) = q^{\pm H/2}\otimes q^{\pm H/2},\
\Delta(X^\pm) =X^\pm\otimes q^{- H/2}+q^{H/2}\otimes X^\pm ,\\
\veps(q^{\pm H/2})=1,\quad \veps(X^\pm)=0.
\endgathered
\tag 32
$$
The left action $\xi$ of $\UU$ on $\CC$ has been derived explicitly in
\cite{\Grass} for a more general case (see also \cite{\JMPII}). After
some rescaling it reads
$$
\xi(q^{\pm H/2})\cdot\z^n=q^{\pm n}\z^n,\
\xi(X^+)\cdot\z^n=[n]\,\z^{n+1},\
\xi(X^-)\cdot\z^n=-[n]\,\z^{n-1},\quad \forall n\in\Zp.
\tag 33
$$

Now we introduce a mapping $\vphi:\MM\to\CC$ by
$$
\vphi(q^{\pm H/2})=q^{\mp\sg/2}\,1,\
\vphi(X^+)=-q^{-\sg/2}[\sg]\,\z,\ \vphi(X^-)=0
\tag 34
$$
where $\sg$ is a complex parameter. Next one has to verify the
assumptions of Proposition 3. It is convenient to add to $\RR$ two other
dependent relations, namely $X^\pm q^{-H/2}-q^{\pm1}q^{-H/2}X^+=0$,
getting this way a new set of relations $\RR'$. Naturally, $\RR$ and
$\RR'$ define the same algebra. The advantage of this step is,
however, that $\RR'$ obeys the condition
$\Delta(\RR')\subset\RR'\otimes\FF+\FF\otimes\RR'$,
as one can check by a direct computation. According to the second part
of Proposition 3 it suffices to verify the equality
$\tp(\RR')=0$ rather than the more complicated assumption (27). This is
again a matter of a straightforward computation. Applying the
prescription (11) we arrive at formulae for the new action:
$$
q^{\pm H/2}\cdot\z^n=q^{\mp(\sg-2n)/2}\,\z^n,\
X^+\cdot\z^n=q^{-\sg/2}[n-\sg]\,\z^{n+1},\
X^-\cdot\z^n=-q^{\sg/2}[n]\,\z^{n-1} ,
\tag 35
$$
valid for all $n\in\Zp$.

Observe that for $\sg\in\Zp$ the unit generates a finite-dimensional
submodule, $\UU\cdot 1=\spn_\BC\{1,\z,\dots,\z^\sg\}$. This is how we
get finite-dimensional irreducible representations of the algebra
$\UU$. In fact, this is a consequence of a more general result about
representations of $\Uq(\gtg)$ \cite{\Rosso,\Lusztig}. Actually, the
unit is in this case a cyclic vector and, at the same time, a lowest
weight vector ($X^-\cdot1=0$) and so the submodule $\UU\cdot 1$ is
unambiguously determined, up to isomorphism, by the lowest weight
given by $q^H\cdot 1=q^{-\sg}\,1$.

 \medpagebreak
 \flushpar {\it 3.2 The general case}
 \smallpagebreak

In the general case we prefer the description of $\UU=\Uq(\gtg)$ due to
Faddeev--Reshetikhin--Takhtajan \cite{\FRT}. The generators are arranged
in respectively upper and lower triangular matrices $L^+$
and $L^-$ of size $N\times N$ and obeying the defining relations
$$
\gathered
R_{12}\Lpm_2\Lpm_1=\Lpm_1\Lpm_2R_{12},\quad
R_{12}\Lp_2\Lm_1=\Lm_1\Lp_2R_{12},\\
\diag(\Lp)\diag(\Lm)=\diag(\Lm)\diag(\Lp)=\BI,\quad
\det(\Lp)=1.
\endgathered
\tag 36
$$
Furthermore, the comultiplication and the counit are determined by
$$
\Delta(\Lpm)=\Lpm\dot{\otimes}\Lpm,\quad
\veps(\Lpm)=\BI
\tag 37
$$
(as usual, $(A\dot{\otimes}B)_{ij}:=\sum_k A_{ik}\otimes B_{kj}$).

Here $R$ is the standard R-matrix obeying the Yang-Baxter equation
$R_{12}R_{13}R_{23}=\mathbreak R_{23}R_{13}R_{12}$
(c.f. \cite{\JimboR} and also \cite{\DrinfeldP, \FRT}) and $\det(\Lp)$
is just the product of diagonal entries. Let us recall that $R$ is lower
triangular (in the lexicographic ordering of indices),
$R_{12}^{\,t}=R_{21}$, and $R_{jk,jt}=0$ for $k\neq t$.

For the series $B_\ell$,
$C_\ell$ and $D_\ell$ there are two additional relations, namely
$(C\Lpm C^{-1})^t=(\Lpm)^{-1}$ where $C$ is a $q$-deformation of the
"classical" matrix $C_0$ (occurring, for example, in the formula (2)).
After having introduced the $N^2\times N^2$ matrix $K$ defined by
$ K_{jk,st}:=C^t_{jk}C^{-1}_{\,st}$
one can rewrite these relations as
$$
K_{12}\Lpm_2\Lpm_1 = K_{12}.
\tag 38
$$
The matrix $K$ is related to the R-matrix by the equality
$$
R_{12}-R_{21}^{-1}=(q-q^{-1})(P-K_{12})
\tag 39
$$
where $P$ stands for the flip operator
($P_{jk,st}=\delta_{jt}\delta_{ks}$).
Taking (39) for the definition of $K$ one finds that $K=0$ for the
series $A_\ell$. Naturally, the conditions (38) become in this case
trivial.

There exist several useful identities involving the matrices $K$ and
$R$ \cite{\JAlg}. Here we mention just the equalities
 $$
 K_{12}R_{31}^{-1}=K_{12}R_{32},\quad
 K_{12}R_{23}^{-1}=K_{12}R_{13},
 \tag 40
 $$
and the implication
 $$
 K_{12}D_1D_2=K_{12} \Longrightarrow
 R_{12}D_1D_2=D_1D_2 R_{12}
 \tag 41
 $$
valid for any complex diagonal matrix $D$.
Particularly, set
 $$
 Q := \diag(R) \quad(\text{then}\ Q_{12}=Q_{21}).
 \tag 42
 $$
 It holds true that
$$
K_{12}Q_{13}Q_{23}=K_{12}\quad\text{and}\quad
R_{12}Q_{13}Q_{23}=Q_{13}Q_{23}R_{12}.
\tag 43
$$

Let us now describe the quantized big cell for the generic coadjoint
orbit of $K$
regarded as a complex manifold \cite{\JSII}. The generators
$z_{jk}$, $1\le j<k\le N$, can be arranged in an upper-triangular
matrix $Z$ as given in (1). The commutation relations then read
$$
R_{12}Q^{-1}Z_1QZ_2 = Q^{-1}Z_2QZ_1R_{12}.
\tag 44
$$
Apart of this one imposes an additional "orthogonality" condition
(trivial for the series $A_\ell$):
$$
K_{12}Q^{-1}Z_1QZ_2 = K_{12}
\tag 45
$$
(this is a simplified but equivalent form to that given in the formula
(7.14) in \cite{\JSII} and in the formula (3.10) in \cite{\JAlg}).
But as one can check by a simple manipulation (45) is already a
consequence of (44) and need not be accounted.
For the series $A_\ell$, $C_\ell$ and $D_\ell$ the matrix
$Q=\diag(R)$ commutes with $R$ and so (44) can be simplified to
$RZ_1QZ_2=Z_2QZ_1R$. However this is not the case for the series
$B_\ell$ (this fact was not recognized in \cite{\JSII}) and so one has
to keep the general form (44). As already mentioned, here we construct
the algebra $\CC$ as being generated by the "antiholomorphic"
generators $z^\ast_{jk}$ arranged in the matrix $Z^\ast$
-- the Hermitian adjoint to $Z$. The corresponding commutation
relation is the Hermitian adjoint to (44), namely
$$
R_{12}\Z_2Q\Z_1Q^{-1} = \Z_1Q\Z_2Q^{-1}R_{12}.
\tag 46
$$

The left action $\xi$ is dual to the right quantum dressing
transformation $\gR:\CC\to\mathbreak\CC\otimes\Funq(G)$.
Here $\Funq(G)$ is the
Hopf algebra of quantum functions living on the group $G$
and it is generated by entries of a matrix $T$ -- the vector
corepresentation of $\Funq(G)$. The
dressing transformation of the holomorphic part formally coincides with
the classical action (3), namely $\gR(Z)=(ZT)_{(+)}$ where
on the right hand side
we have identified $\CC$ with $\CC\otimes1$ and $\Funq(G)$ with
$1\otimes\Funq(G)$. To get the dressing transformation of $\Z$ one can
simply apply the $\ast$-involution. But before doing it one has to
pass
from $\Funq(G)$ to the compact form $\Funq(K)$ which means nothing but
introducing a $\ast$-involution on $\Funq(G)$ by $T^\ast:=T^{-1}$.

The left action $\xi$ is defined by
$$
\xi_x\cdot f := (\id\otimes\langle x,\cdot\rangle)\,\gR(f).
\tag 47
$$
The dual pairing between $\Uq(\gtg)$ and $\Funq(G)$ is prescribed on the
generators as follows \cite{\FRT}:
$$
\langle\Lp_1\dot{,}\,T_2\rangle = R_{21},\quad
\langle\Lm_1\dot{,}\,T_2\rangle = R_{12}^{-1}.
\tag 48
$$
A straightforward computation then gives the desired action:
$$
\xi(\Lp_1)\cdot \Z_2 = R_{21}^{-1}\Z_2Q, \quad
\xi(\Lm_1)\cdot\Z_2 = \Z_1Q\Z_2Q^{-1}(\Z_1)^{-1}.
\tag 49
$$
It can be extended to an arbitrary element from $\CC$ with the aid of
Leibniz rule (7) and the prescription for comultiplication (37).

Let us specify the mapping $\vphi$ on the generators:
$$
\vphi(\Lp)=D^{-1},\quad
\vphi(\Lm)=\Z D^2(\Z)^{-1}D^{-1}
\tag 50
$$
where $D$ is an arbitrary complex diagonal matrix obeying the conditions
$$
\det(D)=1\quad\text{and}\quad K_{12}D_1D_2=K_{12}
\tag 51
$$
(the former one follows from the latter one in the case of the series
$B_\ell$, $C_\ell$ and $D_\ell$). It is easy to show that the set of
defining relations $\RR$ corresponding to the equalities (36) and (38)
obeys the condition
$\Delta(\RR)\subset\RR\otimes\FF+\FF\otimes\RR$.
Thus one can again apply the second part of Proposition 3
to conclude that it suffices to verify the equality
$\tp(\RR)=0$ rather than the assumption (27). This is
a matter of a straightforward computation
(based on the rule (23)) to find that
$$
\gathered
\tp(\Lp_1\Lp_2)=D_1^{-1}D_2^{-1},\quad
\tp(\Lm_1\Lm_2)=\Z_1Q\Z_2D_1^2D_2^2
(\Z_2)^{-1}Q^{-1}(\Z_1)^{-1}D_1^{-1}D_2^{-1},\\
\tp(\Lp_2\Lm_1)=
R_{12}^{-1}\Z_1D_1^2(\Z_1)^{-1}D_1^{-1}D_2^{-1}R_{12},\quad
\tp(\Lm_1\Lp_2)=\Z_1D_1^2(\Z_1)^{-1}D_1^{-1}D_2^{-1}.
\endgathered
\tag 52
$$
The equality $\tp(\RR)=0$ follows immediately from (52) and from the
properties of matrices $R$ and $K$ as mentioned above.

Consequently we conclude from Proposition 1 that there exists a new left
action of $\UU$ on $\CC$, $x\otimes f\mapsto x\cdot f$, for which the
diagonal matrix $D$ plays the role of a parameter. Moreover we know that
$x\cdot1=\vphi(x)$ for all $x\in\UU$; particularly this concerns the
entries of $\Lp$ and $\Lm$ (c.f. (50)).

The last observation is devoted to the cyclic submodule $\UU\cdot1$
generated by the unit. Let us recall the structure of the matrices $\Lp$
and $\Lm$ \cite{\FRT}.
The diagonal entries have the form $q^H$ where $H$ is an
element from the Cartan subalgebra $\gh\subset\gtg$ while the entries
above the diagonal of $\Lp$ are proportional to the negative root
vectors $X^-_\alpha$ (or their $q$-commutators) and the entries below
the diagonal of $\Lm$ are proportional to the positive root vectors
$X^+_\alpha$. Hence the unit is a lowest weight vector
($X^-_\alpha\cdot1=0$ for all simple roots $\alpha$) and the
corresponding lowest
weight is determined by the matrix $D$ in accordance with the equality
$\diag(\Lp)\cdot1=D^{-1}$. This is why we can refer in this case, too,
to the general result \cite{\Rosso, \Lusztig} according to which the
cyclic submodule $\UU\cdot1$ is unambiguously determined by the lowest
weight. This implies that for a discrete set of matrices $D$
corresponding to lowest weights
$-(n_1\omega_1+\dots+n_\ell\omega_\ell)$,
with $n_i\in\BZ_+$ and
$\{\omega_1,\dots,\omega_\ell\}\subset\gh^\ast$ being the set of
fundamental weights, the submodule $\UU\cdot1$ is finite-dimensional and
irreducible.

Let us also note that the particular case of $\gtg=\goth{so}(5)$ has been
treated as an example in \cite{\CJP} with the computations carried out
up to the end.

 \medpagebreak
 \flushpar {\it 3.3 Twisted adjoint action}
 \smallpagebreak

Here we wish to give another description of the preceding example while
abandoning the geometric terminology and relying instead on the notion
of a Verma module. The construction of the modified action presented in
Section 2 then yields exactly the so called twisted adjoint action given
in the book \cite{\Joseph}, \S5.3.10 (this fact has been pointed out to the
author by a referee). We warn the reader however that, if compared with
\cite{\Joseph}, the role of the subalgebras $\gb_+,\gb_-\subset\gtg$
is interchanged and the generators of $\gtg$ are partially rescaled. The
description below is rather brief and with some details omitted.

For the generators of $\Uq(\gtg)$ we chose $e_i=q^{H_i/2}\,X_i^+$,
$f_i=X_i^-\,q^{-H_i/2}$, $t_i^{\pm1}=q^{\pm H_i}$, with the index $i$
enumerating a set of simple roots $\{\alpha_i\}_i$. Thus the defining
relations read
$$
[\,e_i,f_j\,]=\delta_{ij}\,\frac{t_i-t_i^{\;-1}}{q-q^{-1}},\
t_i e_j t_i^{\;-1}=q^{\langle\alpha_i,\alpha_j\rangle}\,e_j,\
t_i f_j t_i^{\;-1}=q^{-\langle\alpha_i,\alpha_j\rangle}\,f_j ,
\tag 53
$$
plus the quantum Serre relations. Let us also recall the comultiplication,
$$
\Delta(e_i)=e_i\otimes1+t_i\otimes e_i,\
\Delta(f_i)=f_i\otimes t_i^{\;-1}+1\otimes e_i,\
\Delta(t_i)=t_i\otimes t_i  .
\tag 54
$$
In this particular case we shall need the antipode which is given by
$$
\sigma(e_i)= - t_i^{\;-1}e_i,\  \sigma(f_i)= - f_i t_i,\
\sigma(t_i)=  t_i^{\;-1}.
\tag 55
$$
The Hopf subalgebra $\Uq(\gb)\subset\Uq(\gtg)$ is generated by
the elements $e_i,t_i^{\;\pm1}$, and the symbol $\Uq(\gn)$ designates
the subalgebra of $\Uq(\gtg)$ generated  by the elements $e_i$
(no comultiplication is defined).

The counit is given as usual ($\varepsilon(e_i)=\varepsilon(f_i)=0$,
$\varepsilon(t_i)=1$) and its restriction determines a one-dimensional
$\Uq(\gb)$ module denoted by $V_\varepsilon$. Verma module $M(0)$
with highest weight zero is introduced by
$$
M(0)=V_\varepsilon\otimes_{\Uq(\gb)}\Uq(\gtg) .
\tag 56
$$
Since $\varepsilon(g\lo)\varepsilon(g\lt)=\varepsilon(g)$,
$\forall g\in\Uq(\gb)$, $M(0)$ is endowed with the structure of a coalgebra
according to the rule
$$
\Delta(1\otimes f)=(1\otimes f\lo)\otimes(1\otimes f\lt) .
$$
Thus the dual space $M(0)^\ast$ is a unital algebra, with the unit being
induced by the counit in $\Uq(\gtg)$. Furthermore, the right $\Uq(\gtg)$
action on $M(0)$ induces a left $\Uq(\gtg)$ action $\xi$ on $M(0)^\ast$. It is
easy to see that $\xi$ obeys (6), (7). In fact, to avoid ill defined expressions
one considers the subalgebra $M(0)^\ast_f\subset M(0)^\ast$ formed by elements
$w$ with the property $\dim\xi(\Uq(\gh))\cdot w<\infty$ (the action of the Cartan
subalgebra $\Uq(\gh)\subset\Uq(\gtg)$ is
required to be locally finite). The $\Uq(\gtg)$-module
algebra $M(0)^\ast_f$ is nothing but the algebra $\CC$ used in the previous
subsection.

According to Corollary 5.3.6 of \cite{\Joseph} the $\Uq(\gtg)$ module
$M(0)^\ast_f$ is isomorphic to the algebra $\Uq(\gn)$. The action on the
latter module is induced by the adjoint action
$$
\ad_x y=x\lo y\, \sigma(x\lt),\ \forall x,y\in\Uq(\gtg) .
\tag 57
$$
In more detail, consider the filtration $F$ on $\Uq(\gtg)$ given by
$\deg(f_i)=1$, $\deg(e_i)=0$ and $\deg(t_i)= - 1$. The filtration is ad-invariant
and consequently there is an induced action of $\Uq(\gtg)$ on
$\gr_F\Uq(\gtg)$. One observes that $\gr_F\Uq(\gn)$ is ad-invariant and can be
identified with $\Uq(\gn)$ as an algebra. It is not difficult to find that the
$\Uq(\gtg)$ action is prescribed on the generators of $\Uq(\gn)$ as follows:
$$
\xi(e_i)\cdot e_j=e_i e_j-q^{\langle\alpha_i,\alpha_j\rangle}\,e_j e_i,\
\xi(f_i)\cdot e_j=\frac{\delta_{ij}}{q-q^{-1}}\,1,\
\xi(t_i)\cdot e_j=q^{\langle\alpha_i,\alpha_j\rangle}\,e_j
\tag 58
$$
(and $\xi(x)\cdot1=\varepsilon(x)\,1$, $\forall x\in\Uq(\gtg)$). The action extends
to the whole algebra $\Uq(\gn)$ with the aid of Leibniz rule (7).

Finally, as described in Proposition 1, the action $\xi$ admits a modification
with the aid of a mapping $\varphi:\Uq(\gtg)\to\Uq(\gn)$. The mapping is
unambiguously defined by its values on generators:
$$
\vphi(e_i)=(1-q^{2\langle\lambda,\alpha_i\rangle})\,e_i,\
\vphi(f_i)=0,\ \vphi(t_i)=q^{\langle\lambda,\alpha_i\rangle}\,1
\tag 59
$$
where $\lambda\in\gh^\ast$ is a weight.


\vskip 12pt
 \noindent{\bf Acknowledgements.}
 The author is indebted to the referee for his comments.
 Partial support from the grant
 202/96/0218 of Czech Grant Agency is gratefully acknowledged.

 \enddocument